\documentclass[letterpaper, 10 pt, conference]{ieeeconf}
\IEEEoverridecommandlockouts
\usepackage{cite}
\usepackage{amsmath,amssymb,amsfonts}
\usepackage{algorithmic}
\usepackage{graphicx}
\usepackage{textcomp}
\usepackage{xcolor}
\usepackage{diagbox}
\usepackage{graphicx} 
\usepackage{epstopdf}
\usepackage{algorithmic}
\usepackage{algorithm}

\def\BibTeX{{\rm B\kern-.05em{\sc i\kern-.025em b}\kern-.08em
    T\kern-.1667em\lower.7ex\hbox{E}\kern-.125emX}}

\makeatletter
\newcommand{\removelatexerror}{\let\@latex@error\@gobble}
\makeatother

\title{\LARGE \bf
Dealing with collinearity in large-scale\\ linear system identification using Bayesian regularization\\
}

\author{Wenqi Cao$^{1}$ and Gianluigi Pillonetto$^{2}$
\thanks{$^{1}$Department of Automation, Shanghai Jiao Tong University, Shanghai, China. {\tt\small wenqicao@sjtu.edu.cn}}%
\thanks{$^{2}$Department of Information Engineering, University of Padova, Padova, Italy. {\tt\small giapi@dei.unipd.it}}%
}

\begin{document}
\maketitle

\begin{abstract}
We consider the identification of large-scale linear and stable dynamic systems
whose outputs may be the result of many correlated inputs.
Hence, severe ill-conditioning may affect the estimation performance.
This is a scenario often arising when modeling complex physical systems
given by the interconnection of many sub-units where feedback
and algebraic loops can be encountered.
We develop a strategy based on Bayesian regularization
where any impulse response is modeled as the realization of a zero-mean Gaussian process.
The stable spline covariance is used to include information
on smooth exponential decay of the impulse responses.
We then design a new Markov chain Monte Carlo scheme
that deals with collinearity and is able to efficiently reconstruct the posterior of the impulse responses.
It is based on a variation of Gibbs sampling which updates possibly overlapping blocks of the parameter space
on the basis of the level of collinearity affecting the different inputs.
Numerical experiments are included to test the goodness of the approach
where hundreds of impulse responses form the system and inputs correlation may be very high.
\end{abstract}


\section{Introduction}

Large-scale dynamic systems arise in many scientific fields
like engineering, biomedicine and
neuroscience \cite{Hagmann2008,Hickman2017,Pagani2013,Prando2020}.
Modeling these complex physical systems
is crucial for prediction and control purposes \cite{Ljung:99,Soderstrom}.
They can be often interpreted as networks composed of
a large set of interconnected sub-units.
One can also describe them through many nodes which can communicate each other
through modules driven by measurable inputs or noises \cite{Chiuso2012,Materassi2010,VdH2013,Yue2021,CPL21}.

We will assume that linear dynamics underly our large-scale dynamic system
focus 
on the identification of Multiple Inputs Single Output (MISO) models.
Hence, system impulse responses have to be estimated from input-output data.
These latter are assumed to be the result of many (measurable) inputs
which can be highly correlated and possibly also poorly exciting. 
Such scenario is often encountered when modeling dynamic networks
since modules interconnections give rise to feedback
and algebraic loops \cite{Bazanella2017,Goncalves2008,Hendrickx2019,Weerts2018,Fonken2020,Ramaswamy2021b}. 
This complicates considerably the identification problem:
low pass and almost collinear inputs lead to severely ill-conditioned estimation problems \cite[Chapter 3]{CollBook}\cite{ChiusoPicci04},
as also described in many engineering applications like \cite{LLB16,MPJR10,PV97}.
In addition, the use of  the classical approach to system identification \cite{Ljung:99,Soderstrom}, where
different parametric structures have to be postulated, 
raises several difficulties. For instance, the use of rational transfer functions
to describing each impulse response leads to high-dimensional nonconvex optimization problems.
In addition, the system dimension has typically to be estimated from data,
e.g. using AIC or cross validation \cite{Akaike:74,Hastie01},
and this further complicates the problem
due to the combinatorial nature of model order selection.

The problem of collinearity in linear regression has been discussed for several decades \cite{CollBook,Srisa21,Zhang18}.
Two types of methods are mainly considered. One consists of eliminating some related regressors, the other of exploiting
special regression methods like LASSO, ridge regression or partial least squares \cite{Srisa21}.
The main idea underlying these approaches is to eliminate or decrease the influence of correlated variables.
This strategy is useful when the goal is to obtain a model useful for predicting future data
from inputs with statistics
similar to those in the training set.
In this paper, we instead follow a different route based on
Bayesian identification of dynamic systems via Gaussian regression \cite{Rasmussen,Bell2004,SurveyKBsysid,BottegalQuant2017}.
In particular, we will couple the stable spline prior
proposed in \cite{SS2010,PillACC2010,SS2011} 
with a Markov chain Monte Carlo (MCMC) strategy \cite{Gilks} suited to overcome collinearity.
In particular, the stable spline kernel is a nonparametric description of stable systems
which has shown some important advantages in comparison
with classical parametric prediction error methods \cite{Ljung:99,SurveyKBsysid,PillonettoBook2022}.
It includes information on
smooth exponential decay and depends on two different kinds of hyperparameters:
a scale factor and a common decay rate \cite{PillonettoTAC2021,PillonettoMLrob2015,ARAVKIN2012}.
The main novelty in this paper is that the stable spline kernel will be implemented
using a Full Bayes approach where the posterior of impulse responses
is sampled through an MCMC strategy able to deal with collinearity.
This is obtained by designing a variation of Gibbs sampling
able to visit more frequently those parts of the parameter space 
more correlated.

The structure of the paper is as follows. Section II reports the problem statement by introducing the measurement model and the prior model on the unknown impulse responses. In Section III, we describe the system identification procedure, discuss the issues regarding the selection of the overlapping blocks to be updated during the MCMC simulation and the impact of stable spline scale factors on the estimation performance. Section IV illustrates two numerical examples. Conclusions then end the paper.

\section{Problem Formulation}\label{secProb}
The measurement model is
\begin{equation}\label{MM1}
y_i =  \sum_{k=1}^m F_k u_k + e_i, \quad i=1,\ldots,n
\end{equation}
where the $F_k$ are stable rational transfer functions sharing a common denominator,
$y_i$ is the output measured at $t_i$ while $u_k$ is the known input entering the $k$-th channel.
Finally, the random variables $e_i$ form a white Gaussian noise
of variance $\sigma^2$. The problem is to estimate the  $m$
transfer functions from the input-output samples. The difficulty we want to face is that
the number of unknown parameters can be relatively large w.r.t. the data set size $n$ 
and some of the $u_k$
may be highly correlated.

We rewrite
\eqref{MM1} in the matrix-vector form by assuming that
FIR models of (possibly large) order $p$ can well approximate each $F_k$.
In particular, $G_k \in \mathbb{R}^{n \times p}$ are suitable Toeplitz matrices containing
past input values such that
\begin{equation}\label{MM2}
Y = \left( \sum_{k=1}^m G_k \theta_k \right) + E = G \theta + E
\end{equation}
where the column vectors $Y,\theta_k$ and $E$ contain, respectively,
the $n$ output measurements, the $p$ impulse response coefficients defining the
$k$-th impulse response and the $n$ i.i.d.
Gaussian noises. Finally, on the rhs $G=[G_1 \ldots G_m]$
while $\theta$ is the
column vector which contains all the $\theta_k$.

Following the Bayesian approach to linear system identification documented in
\cite{SS2010,SurveyKBsysid}, for a known covariance, the $\theta_k$ are seen as Gaussian vectors
with zero-mean and covariances proportional to the stable spline matrix  $K \in \mathbb{R}^{p \times p}$
that encodes smooth exponential decay information.
In particular, the $i,j$-entry of $K$ is
\begin{equation}\label{SS}
K(i,j) = \alpha^{\max(i,j)}. 
\end{equation}
where $\alpha$ regulates how fast the impulse responses are expected to go to zero.
This parameter will be assumed to be known in what follows
to simplify the exposition.

Differently from \cite{SS2010}, the covariances of the $\theta_k$ depend on scale factors $\lambda_k$
which are seen as independent random variables following an improper Jeffrey's distribution
\cite{Jeffreys1946}. Such prior in practice includes only nonnegativity information and is specified by
\begin{equation}\label{Hp1}
\mathrm{p}(\lambda_k) \sim \frac{1}{\lambda_k}
\end{equation}
where here, and in what follows, $\mathrm{p}(\cdot)$ denotes a probability density function.
Later on, we will also see that constraining all the scale factors to be the same can be crucial to deal with collinearity.
The noise variance $\sigma^2$ is also a random variable (independent of $\lambda_k$) following the Jeffrey's prior.
Finally, for known $\lambda_k$ and $\sigma^2$, all the $\theta_k$ and the noises in $E$ are assumed mutually independent.

\section{Gibbs Sampling with Overlapping Blocks}\label{secAlg}

\subsection{Bayesian regularization}
From Section~\ref{secProb}, one obtains
\begin{subequations}\label{priDis}
\begin{equation}
    Y \mid (\{\theta_k\},\{\lambda_k\},\sigma^2)\sim  \mathcal{N}(G\theta, \sigma^2I),
\end{equation}
\begin{equation}\label{priDisTheta}
    \theta_k \mid \lambda_k \sim \mathcal{N}(0,\lambda_k K), 
\end{equation}
\begin{equation}
    E \mid \sigma^2 \sim \mathcal{N}(0,\sigma^2 I). 
\end{equation}
\end{subequations}

Using Bayes rule, one has 
\begin{subequations}\label{postDis}
\begin{equation}\label{postDisLambdak}
    \lambda_k \mid (Y, \sigma^2, \{\theta_k\}, \{\lambda_j\}_{j\neq k}) \sim {\mathcal{I}_g}(\frac{p}{2},~ \frac{1}{2}\theta_k'K^{-1}\theta_k),
\end{equation}
\begin{equation}
    \sigma^2 \mid (Y,\{\theta_k\}, \{\lambda_k\})\sim  {\mathcal{I}_g} (\frac{n}{2}, \frac{1}{2}{\Vert Y-G\theta\Vert^2}),
\end{equation}
\begin{equation}
    \theta_k \mid (Y,\sigma^2, \{\lambda_k\}, \{\theta_j\}_{j\neq k}) \sim \mathcal{N}(\hat{\mu}_k, \hat{\Sigma}_k),
\end{equation}
where ${\mathcal{I}_g(\cdot,\cdot)}$ denotes the inverse Gamma distribution,
\begin{equation}
    \hat{\mu}_k= \hat{\Sigma}_k\frac{1}{\sigma^2}G_k'(Y-\sum_{j\neq k}G_j\theta_j),
\end{equation}
\begin{equation}
    \hat{\Sigma}_k=(\lambda_k^{-1}K^{-1}+\frac{1}{\sigma^2}G_k'G_k)^{-1}.
\end{equation}
\end{subequations}
The above equations would already allow to implement a Gibbs sampling approach
where, at any iteration $t$, the following samples are generated
\begin{subequations}
\begin{equation}\label{gibbsLambdak}
\begin{split}
   &\lambda_k^{(t)} \mid (Y, \sigma^{2(t-1)},\{\lambda_j^{(t)}\}_{j=1}^{j=k-1}, \{\lambda_j^{(t-1)}\}_{j=k+1}^{j=m},\\ &
   \{\theta_k^{(t-1)}\} ),  \quad \quad \text{for}~k=1,\cdots,m,
\end{split}
\end{equation}
\begin{equation}\label{gibbsSigma2}
    \sigma^{2(t)} \mid (Y,\{\lambda_k^{(t)}\},\{\theta_k^{(t-1)}\}),\quad\text{for}~k=1,\cdots,m,
\end{equation}
\begin{equation}\label{gibbsThetak}
\begin{split}
    \theta_k^{(t)} \mid (Y,\sigma^{2(t)},\{\lambda_k^{(t)}\}, \{\theta_j^{(t)}\}_{j=1}^{k-1}, \{\theta_j^{(t-1)}\}_{j=k+1}^{m}),&\\
    \text{for}~k=1,\cdots,m.&
\end{split}
\end{equation}
\end{subequations}
However, such classical approach
may have difficulties to handle collinearity, leading to a slow mixing of the chain.
This motivates the development described in the next sections.

\subsection{Overlapping blocks}
In our system identification setting, collinearity is related to the
relationship between two inputs $u_i$ and $u_j$, making some posterior regions more difficult to be explored.
We need a collinearity measure that allows our sampling strategy to
focus more on these parts of the parameter space.
In particular, we still generate samples of the single $\theta_i$
according to the equations described above but we want also
to update larger blocks containing couples of impulse responses
which are significantly correlated a posteriori.
A first important index is the absolute value of the correlation coefficient, i.e.
\begin{equation}\label{cij}
    c_{ij}:=\left| \frac{\text{Cov}(u_i,u_j)}{\sqrt{\text{Var}(u_i)\text{Var}(u_j)}}\right|,
\end{equation}
where $\text{Cov} (\cdot,\cdot)$ and $\text{Var}(\cdot)$ denote the covariance and variance, respectively.
It is now crucial to define a suitable function which maps
each $c_{ij}$  into the probability to update the vector containing both $\theta_i$ and $\theta_j$ at once.

\begin{figure}[htbp]
\centering
\includegraphics[scale=0.6]{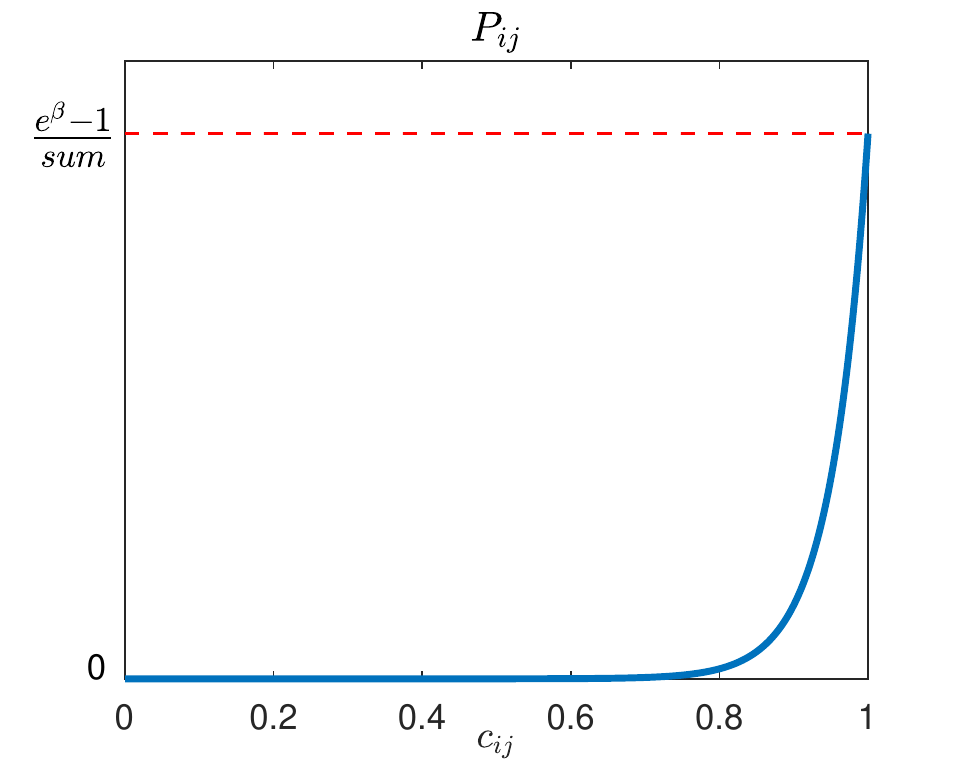}
\caption{$P_{ij}$ with respect to $c_{ij}$}
\label{cp}
\end{figure}

Formally, define $\theta_{ij}:= [\theta_i', \theta_j']'$ and let $P_{ij}$ be the probability of
selecting $\theta_{ij}$ as the block to be updated inside an MCMC iteration.
In particular, we will use an exponential rule
that emphasises  correlation coefficients close to one, i.e.
\begin{equation}\label{Pij}
    P_{ij}=\left\{\begin{array}{cc}
               0,& i=j;  \\
               \frac{e^{\beta c_{ij}}-1}{sum}, &i\neq j;
             \end{array}\right.
\end{equation}
where $sum=\sum_{i<j}(e^{\beta c_{ij}}-1)$, and $\beta$ is a tuning rate.
For instance, the profile of $P_{ij}$ as a function of $c_{ij}$ for $\beta=20$
is displayed in Fig.~\ref{cp}.

Now, letting $G_{ij}=[G_i,G_j]$, the conditional distribution of $\theta_{ij}$ is
\begin{subequations}\label{BlockUpdate}
\begin{equation}
    \theta_{ij}\mid (Y, \sigma^{2},\{\theta_k\}_{k\neq i,j}, \{\lambda_k\}) \sim \mathcal{N}(\hat{\mu}_{ij}, \hat{\Sigma}_{ij})
\end{equation}
where
\begin{equation}
    \hat{\Sigma}_{ij}=(\begin{bmatrix}\lambda_i^{-1}K^{-1}&0\\
    0 & \lambda_j^{-1}K^{-1}\end{bmatrix}
    +\frac{1}{\sigma^2}G_{ij}'G_{ij})^{-1},
\end{equation}
\begin{equation}
      \hat{\mu}_{ij}= \hat{\Sigma}_{ij}\frac{1}{\sigma^2}G_{ij}'(Y-\sum_{k\neq i,j} G_k\theta_k).
\end{equation}
\end{subequations}
In conclusion, 
our scheme updates all the parameters sequentially
by using (\ref{gibbsLambdak}-\ref{gibbsThetak}). But, in addition,
according to $P_{ij}$, some larger impulse responses blocks are selected and updated.
Following \cite[Chapter 1.3]{Gilks}, the acceptance rate of such overlapping blocks $\theta_{ij}$ is always 1.

\subsection{A common scale factor}
The simulation strategy described in the previous section uses a model
where the covariance of each impulse response depends on a
different scale factor $\lambda_k$.
As described also by simulation results in Section~\ref{secExp},
in presence of strong collinearity,
such statistical description may make the chain nearly reducible \cite[Chapter 3]{Gilks}, 
i.e. unable to visit the highly correlated parts of the posterior. 
 The problem can be overcome by reducing the model complexity,
i.e. assigning a common scale factor to all the impulse responses covariances.
Hence, the distribution \eqref{postDisLambdak} becomes
\begin{equation}\label{postDisLambda}
    \lambda \mid (Y, \sigma^2, \{\theta_k\}) \sim {\mathcal{I}_g}(\frac{np}{2},~ \frac{1}{2}\sum \theta_k'K^{-1}\theta_k)
\end{equation}
and the update of $\lambda$, previously define by \eqref{gibbsLambdak}, becomes
\begin{equation}\label{gibbsLambda}
    \lambda^{(t)} \mid (Y, \sigma^{2(t-1)}, \{\theta_k^{(t-1)}\} ).
\end{equation}
Hence,  at step $t$ of the MCMC algorithm, when $\theta_{ij}$ is selected it is updated according to
\begin{equation}\label{gibbsThetaij}
     \theta_{ij}^{(t)} \mid (Y, \sigma^{2(t)},\lambda^{(t)}, \{\theta_k^{(t)}\}_{k\neq i,j}),
\end{equation}
i.e. using \eqref{BlockUpdate} with $\lambda_i=\lambda_j=\lambda^{(t)}$.
We are now in a position to
introduce our MCMC strategy for linear system identification under collinear inputs,
which relies on
Gibbs samplings with overlapping blocks (GSOB). Using $n_{MC}$ to indicate the number of MCMC steps, $n_{OB}$ the number of overlapping blocks which are selected and updated at each iteration, and $n_B$ the number of burn-in \cite[Chapter 7]{Gilks},  such procedure is summarized in Algorithm~1.

\begin{figure}[!t]
  \label{algGSOB}
  \renewcommand{\algorithmicrequire}{\textbf{Input:}}
  \renewcommand{\algorithmicensure}{\textbf{Output:}}
  \removelatexerror
  \begin{algorithm}[H]
    \caption{Gibbs samplings with overlapping blocks using one 
    common scale factor (GSOB)}
    \begin{algorithmic}[1]
      \REQUIRE Measurements $G$, $Y$; initial values $\lambda^{(0)}$, $\sigma^{2(0)}$, $\{\theta_k^{(0)}\}$. $\alpha$, $\beta$.
      \ENSURE Estimate $\hat\theta$.
      \STATE Calculate $\{P_{ij}\}$ from input data for $i\leq j$ in \eqref{Pij};
      \FOR {$t=1:n_{MC}$}
      \STATE Sample \eqref{gibbsLambda}, \eqref{gibbsSigma2} and \eqref{gibbsThetak} in sequence;
        \FOR {$k=1:n_{OB}$}
        \STATE Choose a pair $(i,j)$ from the distribution $\{P_{ij}\}_{i\leq j}$;
        \STATE Sample \eqref{gibbsThetaij} and update the values of $\theta_i^{(t)}$ and $\theta_j^{(t)}$.
        \ENDFOR
      \ENDFOR
     \STATE Calculate $\hat\theta$ as the mean of the impulse responses samples from $t=n_B+1$ to $t=n_{MC}$. 
    \end{algorithmic}
  \end{algorithm}
\end{figure}

\section{Simulation examples}\label{secExp}
\subsection{Example 1}
In the first example we shall use one extreme case with two inputs that are exactly the same.
This will illustrate the difficulties due to collinearity, as well as the importance of adopting only a common scale factor $\lambda$
for all the impulse responses. 

\begin{figure}[htbp]
\centering
\includegraphics[scale=0.6]{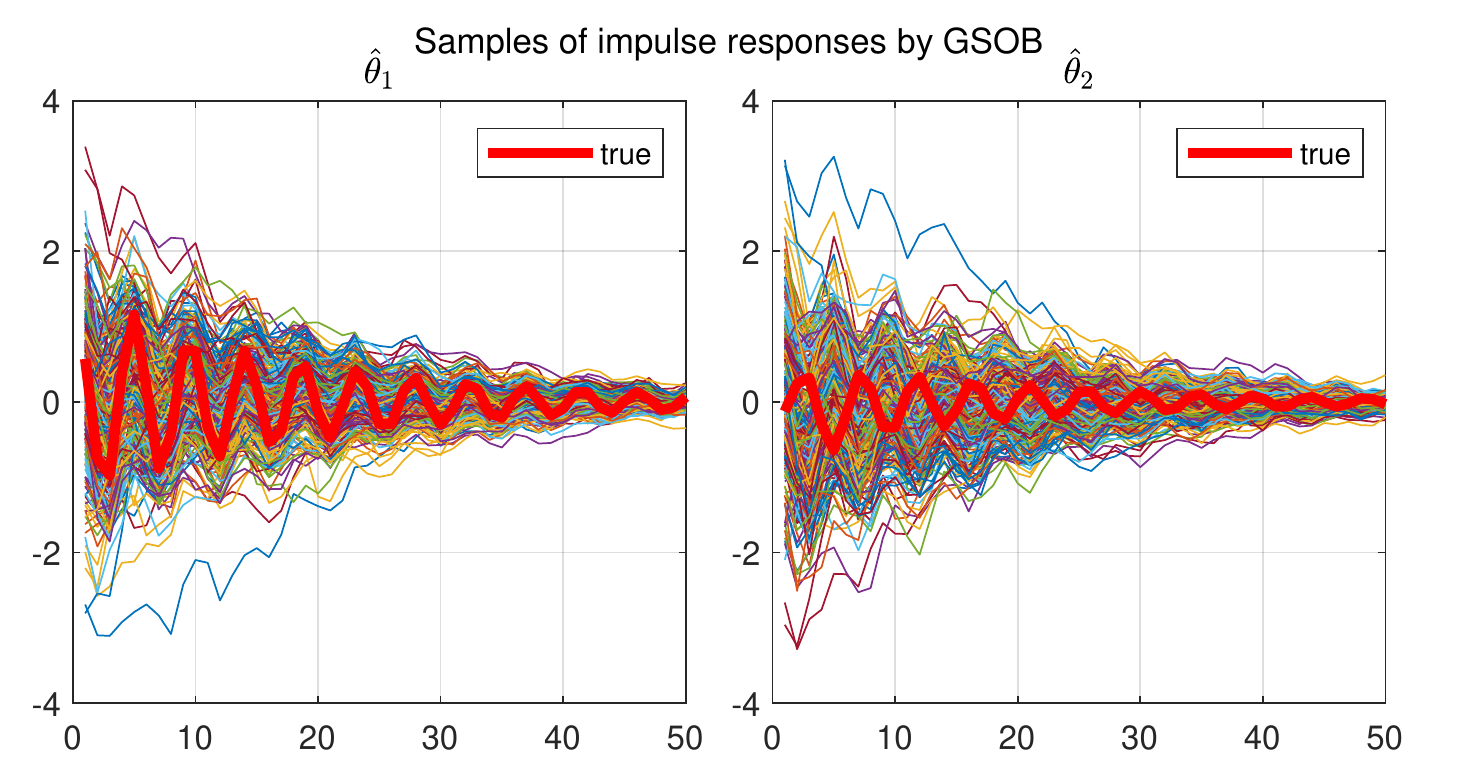}
\caption{Samples of impulse responses after 500 iterations of algorithm GSOB in Example 1.}
\label{ex1_ir_GSOB}
\end{figure}

Let $u_1(t)=u_2(t)$, with the inputs defined by realizations of
white Gaussian noise with $n=500$. We consider two randomly generated transfer functions $F_1(z), F_2(z)$ with a common denominators of degree 5, displayed e.g. in the
two panels of Fig. \ref{ex1_ir_GSOB}.
The measurement noises $e_i$ are independent Gaussian of variance $0.3$.
In the Fisherian framework used by classical system identification this problem is non-identifiable.
But we use model \eqref{MM2} of order $p=50$ and the Bayesian framework
to find impulse responses estimates. In this experiment,
and in Example 2, we set $\alpha=0.9$.

In this example,  $c_{12}=c_{21}=1$, so that
$P_{12}=P_{21}=1$ after normalization. 
We run $n_{MC}=500$ Monte-Carlo iterations in MATLAB to compare the following four algorithms:
\begin{itemize}
  \item GSOB: Gibbs sampling with overlapping blocks using the same scale factor with $n_{OB}=2$;
  \item GSOBd: Gibbs sampling with overlapping blocks using different scale factors $n_{OB}=2$;
  \item GS: Gibbs sampling without overlapping blocks using the same scale factor;
  \item GSd: Gibbs sampling without overlapping blocks using different scale factors.
\end{itemize}

\begin{figure}[htbp]
\centering
\includegraphics[scale=0.6]{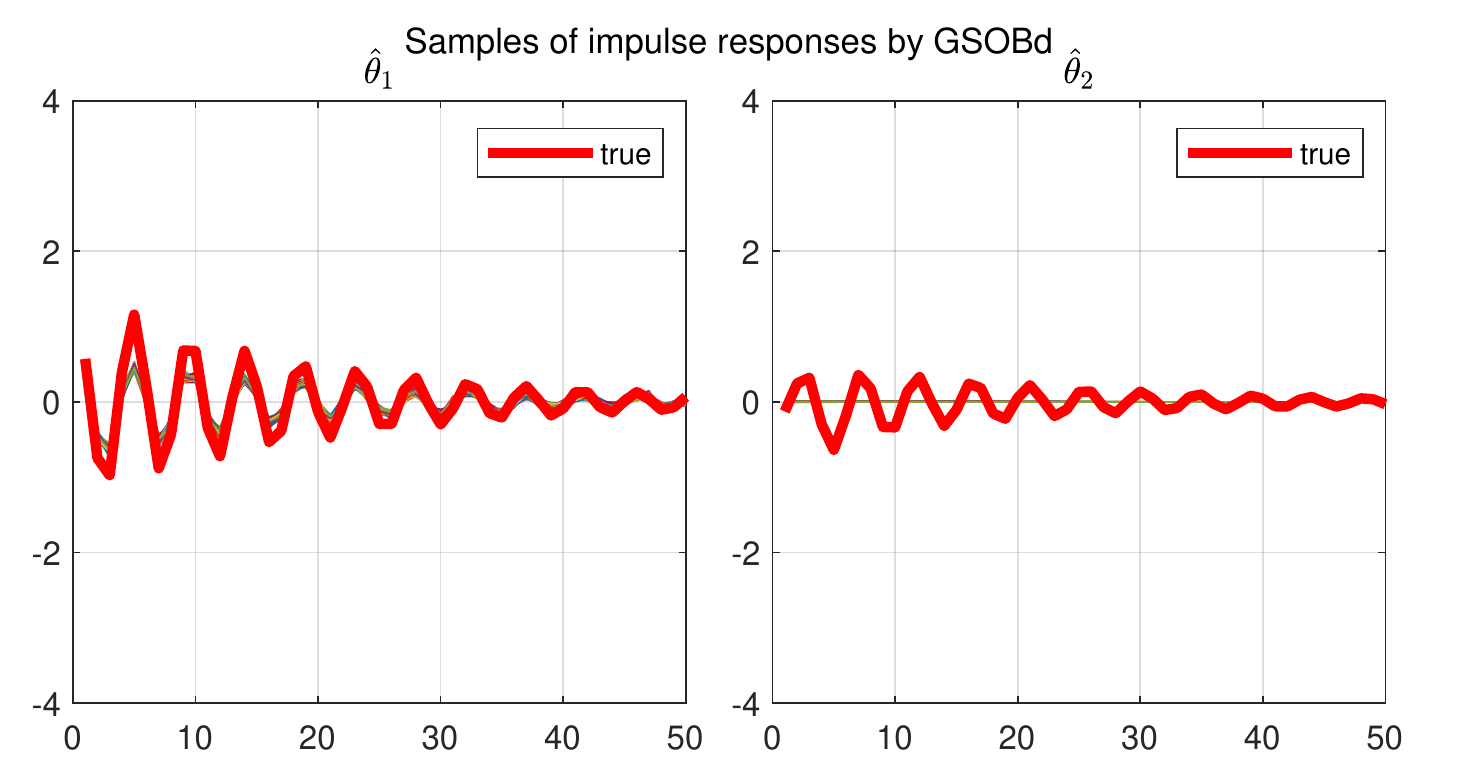}
\caption{Samples of impulse responses after 500 iterations of algorithm GSOBd in Example 1.}
\label{ex1_ir_GSOBd}
\end{figure}

\begin{figure}[htbp]
\centering
\includegraphics[scale=0.6]{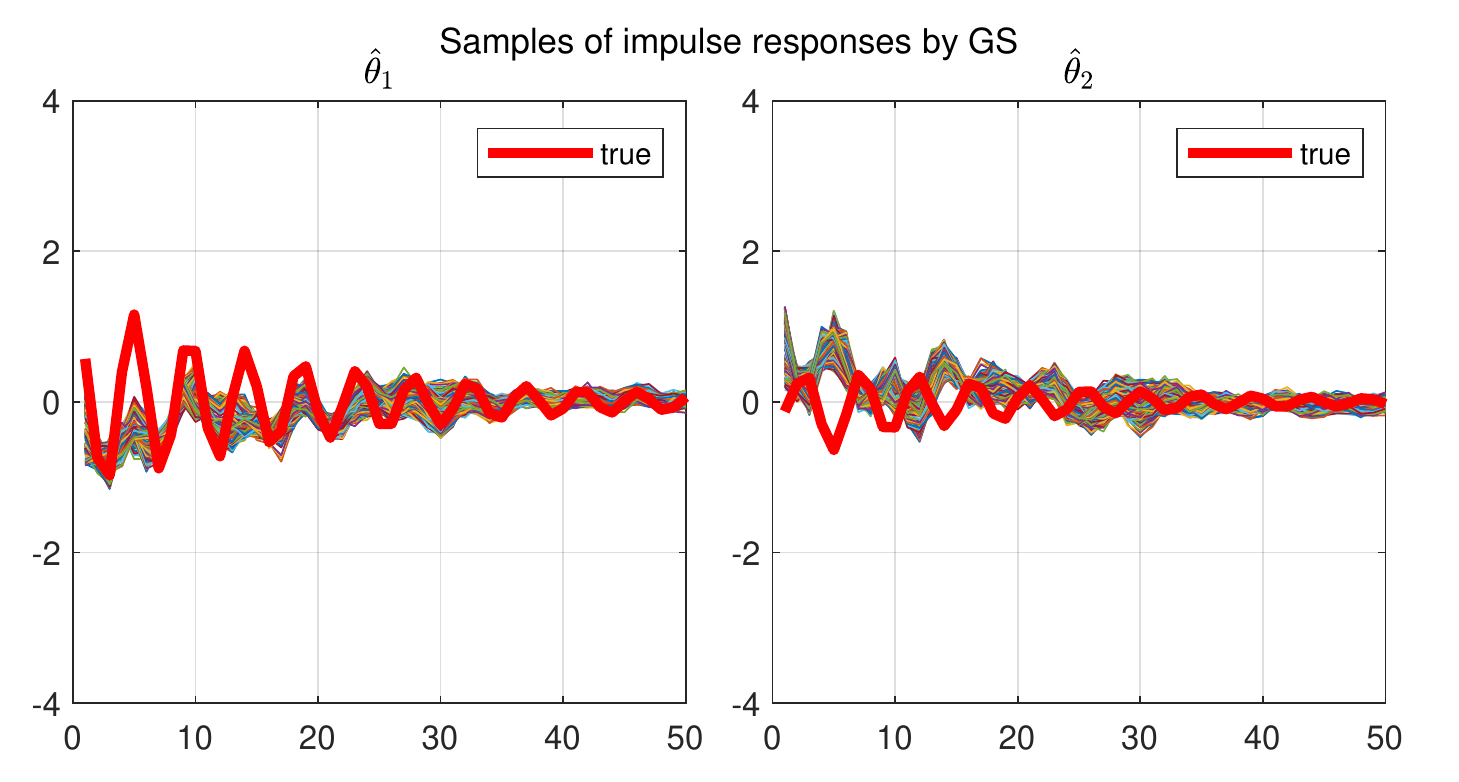}
\caption{Samples of impulse responses after 500 iterations of algorithm GS in Example 1.}
\label{ex1_ir_GS}
\end{figure}

\begin{figure}[htbp]
\centering
\includegraphics[scale=0.6]{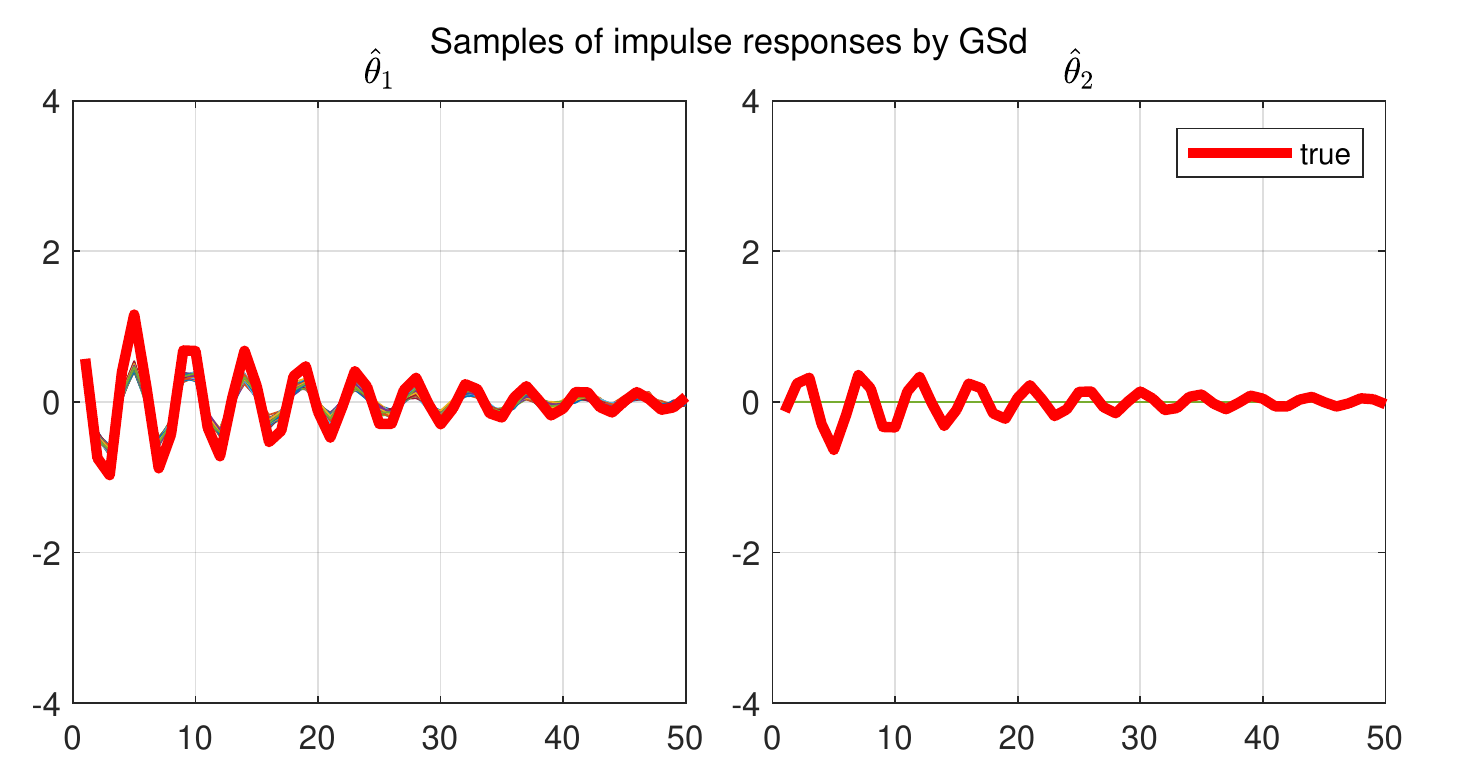}
\caption{Samples of impulse responses after 500 iterations of algorithm GSd in Example 1.}
\label{ex1_ir_GSd}
\end{figure}

The posterior of the two impulse responses is reconstructed 
using samples from MCMC considering the first $50\%$ as burn-in. 
Results coming from the four procedures are displayed in Figs.~\ref{ex1_ir_GSOB}, \ref{ex1_ir_GSOBd}, \ref{ex1_ir_GS} and \ref{ex1_ir_GSd} (true impulse responses are the red thick lines).

After the 500 iterations, only GSOB returns a sampled version of the posterior
able to highlight the non-identifiability issues underlying the problem, see Fig.~\ref{ex1_ir_GSOB}.
Results from classical Gibbs sampling are in Fig.~\ref{ex1_ir_GS}. One can see that samples of $\theta_2$
have not yet reached the same informativeness of GSOB, which may lead to a large prediction error when the input data are independent. The reason is that GS converges much slower than
GSOB under collinearity.

The situation is even more critical when adopting different scale factors.
In fact, Figs.~\ref{ex1_ir_GSOBd} and \ref{ex1_ir_GSd} show that the samples of the second impulse response $\theta_2$ are
all close to zero. The reason is that chains with multiple scale factors may easily become nearly reducible.
At the first MCMC iterations, the algorithm generates an impulse response $\theta_1$ able to explain all the data.
Hence, $\theta_2$ is virtually set to zero and then also $\lambda_2$ becomes very small.
This scale factor becomes a strong (and wrong) prior for $\theta_2$ during the next iterations,
always suggesting the sampler use only $\theta_1$ to fit the output data.

\subsection{Example 2}
In this example we shall use 100 inputs to test the appearance of our algorithm when facing large-scale systems.
This is a situation where dividing the parameter space into many small blocks can have crucial
computational advantages. Otherwise, one should deal with the inversion of very large matrices whose dimension
is the overall number of impulse response coefficients.

We generate the collinear inputs as follows
\begin{equation}\label{col}
 u_i(t)=u_j(t)+r_{ij}(t),
\end{equation}
where $r_{ij}(t)$ is a noise independent of $u_j(t)$.
To increase collinearity between inputs, and the ill-conditioning affecting the problem,
$r_{ij}(t)$ in \eqref{col} is a moving average (MA) process,
\begin{equation}\nonumber
    r_{ij}(t)=v_{ij}(t)-0.8v_{ij}(t-1),
\end{equation}
where $v_{ij}(t) \sim \mathcal{N}(0,\gamma^2_{ij})$.
From \eqref{cij}, one has
\begin{equation}\nonumber
    c_{ij}=\frac{1}{\sqrt{1+\gamma^2_{ij}/(1-0.8^2)}},
\end{equation}
when $u_j(t)$ is a zero mean process with variance 1. 
The level of collinearity can then be tuned by $\gamma_{ij}$.

The data set size is set to $n=10^5$. We introduce correlation
to  $10\%$ of the inputs by letting 
\begin{equation}
    u_{i+1}(t)=u_{i}(t) + r_{(i+1)i}(t),
\end{equation}
for $i=1,\cdots,9$, where the first input $u_1(t)$ is white and Gaussian of variance one, $r_{(i+1)i}(t)$ is obtained by setting the $c_{i(i+1)}=0.99$ for $i=1,2,\cdots,9$. The other $90\%$  of the inputs, i.e. $\{u_i\}_{i=11}^{100}$ are zero-mean standard i.i.d. Gaussian processes.

\begin{figure}[htbp]
\centering
\includegraphics[scale=0.6]{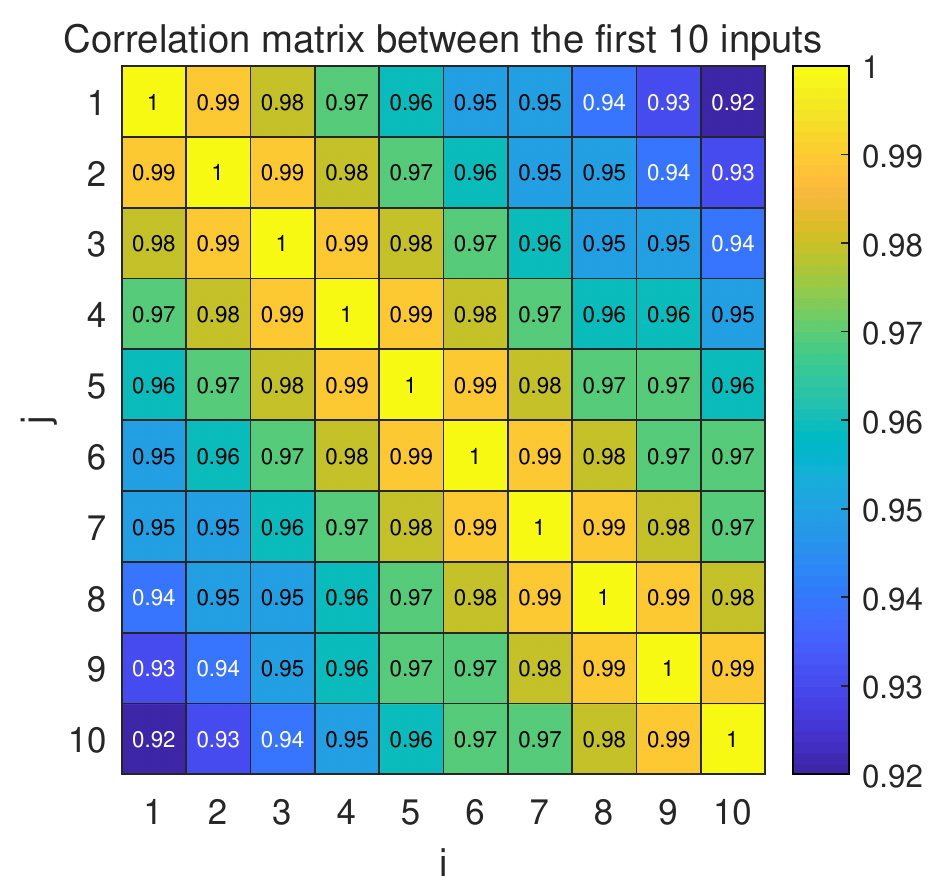}
\caption{$\{c_{ij}\}$ of $\{(u_i,u_j)\}_{i,j=1}^{10}$ in Example 2.}
\label{heatmapCM}
\end{figure}

The different collinearity levels of the first 10 inputs are summarized by the correlation coefficient matrix (approximately calculated from the input realizations) reported in Fig.~\ref{heatmapCM}.
One can see that the input pairs have correlation coefficients ranging from 0.99 to 0.92.

\begin{figure}[htbp]
\centering
\includegraphics[scale=0.6]{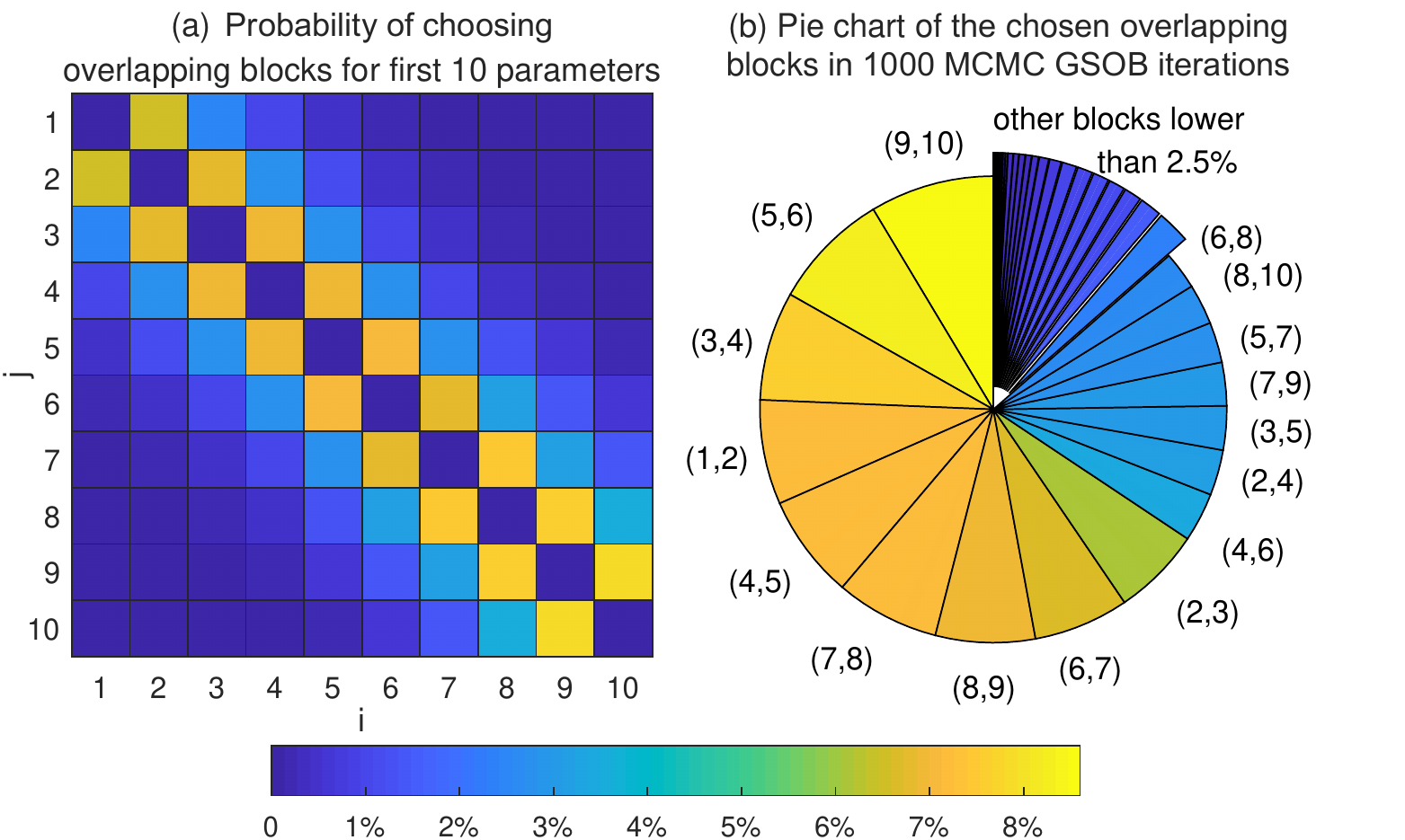}
\caption{(a) Probability $\{P_{ij}\}$ of choosing couples $\{(i,j)\}_{i,j=1}^{10}$ in example 2 of GSOB; (b) pie chart of the selected frequencies of overlapping blocks in 1000 Gibbs iterations of GSOB in Example 2.}
\label{fig6}
\end{figure}

We set $\beta=100$ and then we use \eqref{Pij} to compute the probabilities of selecting the overlapping blocks.
$P_{ij}$ for the first 10 impulse responses are shown in Fig.~\ref{fig6}~(a) with a sharing colorbar at the bottom
($P_{ij}$ for  $i,j>10$ are obviously all close to zero).
In all pairs considered for constructing overlapping blocks $(i,j) (i\neq j)$,
The couples $(i,i+1)$ ($i=1,...,9$) of highest collinearity are assigned almost $7\%$ of the total amount of probability.
On the other hand, e.g. the couple $(1,10)$, which has a correlation coefficient $0.9205$, is given $0.0063\%$.

We generate random transfer functions $F_i(z)~(i=1,\cdots, 100)$ of degree 5 with a common denominator.
Measurement Gaussian noises $e_i$ have variance $0.3$. We adopt model \eqref{MM2}
with each impulse response of order $p=50$ and show results coming only from GSOB and GS,
letting $\alpha=0.9$, $n_{OB}=10$ and $n_{MC}=1000$. 
The pie chart of the selected overlapping blocks' index pairs 
is shown in Fig.~\ref{fig6} (b), where the frequencies of selection of the overlapping blocks are shown in different colors.
As expected, they turn out close to the probabilities $P_{ij}$.
The pair $(i,j)$ will be sampled more frequently along with the growth of the correlation.

\begin{figure*}[htbp]
\centering
\includegraphics[scale=0.8]{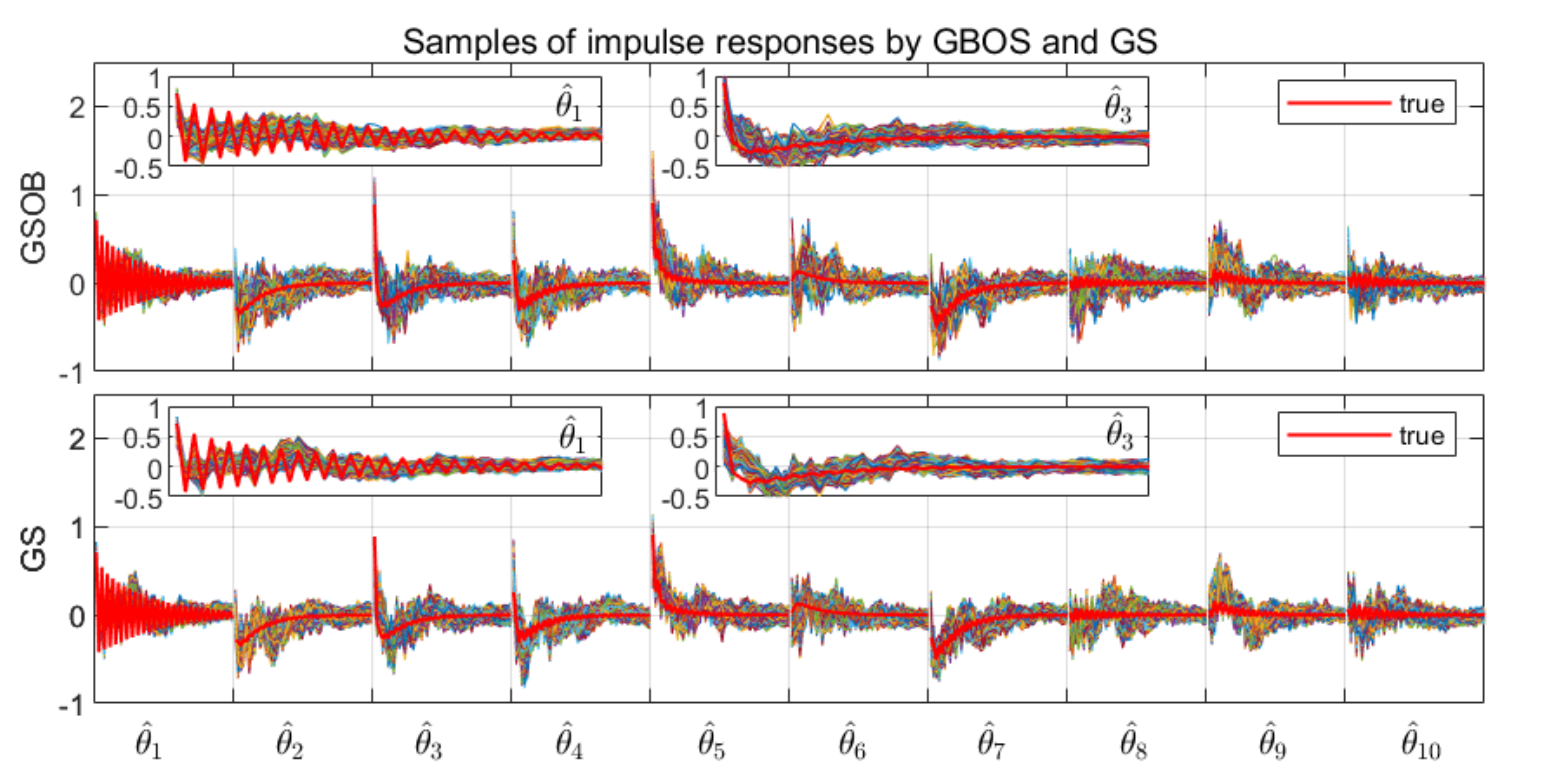}
\caption{Sampled impulse responses of $\{\theta\}_{i=1}^{10}$ after 1000 iterations of algorithms GSOB and GS in Example 2.}
\label{ex2_ir_f10}
\end{figure*}

\begin{figure*}[htbp]
\centering
\includegraphics[scale=0.8]{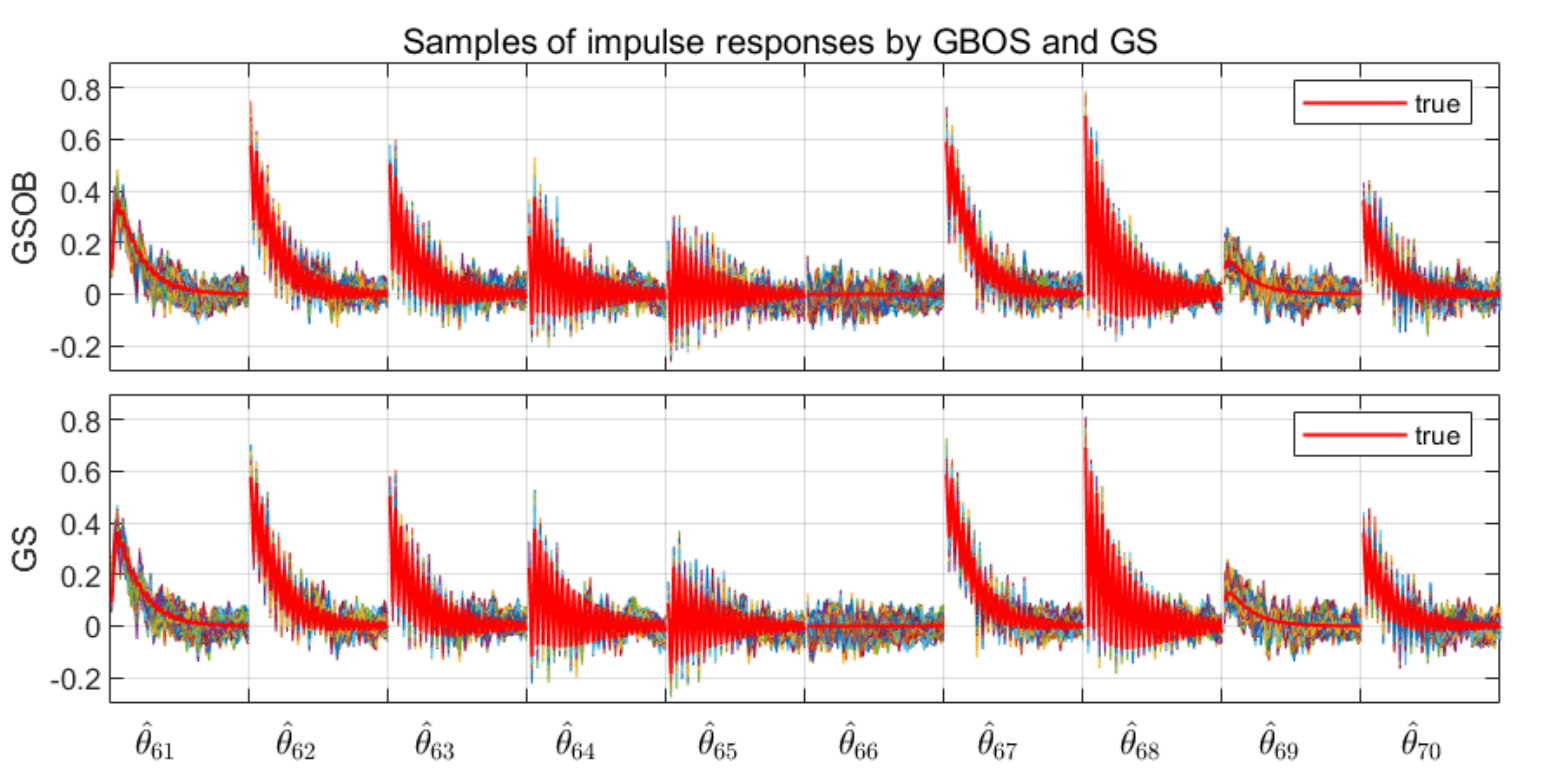}
\caption{Sampled impulse responses of $\{\theta\}_{i=61}^{70}$ after 1000 iterations of algorithms GSOB and GS in Example 2.}
\label{ex2_ir_nc}
\end{figure*}

The last $50\%$ of MCMC realizations are used to reconstruct the
impulse responses posterior in sampled form, as done in Example 1.
Results concerning the first 10 impulse responses (related to the strongly collinear inputs) are shown in Fig.~\ref{ex2_ir_f10}.
One can see that  GSOB works better than GS. For readers' convenience, samples from the posterior of $\theta_1$ and $\theta_3$ are
also zoomed where big differences of the samplings in recovering the true values are shown.
Regarding the other impulse responses (related to inputs which are mutually independent)
the performance of both the algorithms is similar, e.g. samples of $\{\theta_i\}_{i=61}^{70}$
are visible in Fig.~\ref{ex2_ir_nc}.

\begin{figure}[htbp]
\centering
\includegraphics[scale=0.6]{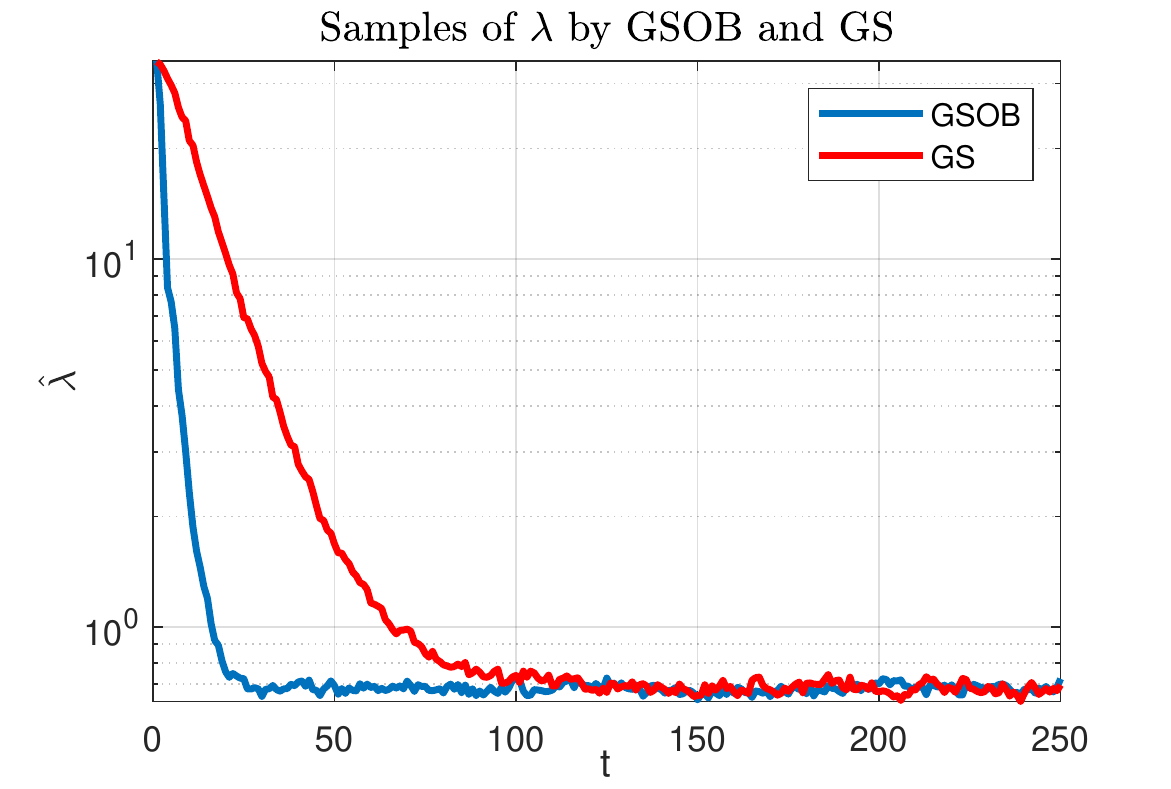}
\caption{Samples of $\lambda$ in the first 250 iterations by GSOB and GS in Example 2.}
\label{ex2_lambda}
\end{figure}

Another fundamental feature of GSOB is its better convergence rate w.r.t. GS.
This can be appreciated e.g. by analyzing the samples from the posterior of  the scale factor
$\lambda$ generated by these two procedures.
They are shown in Fig.~\ref{ex2_lambda}, as a function of the first 250 iterations.
Chain's mixing of GSOB is significantly faster than GS.

\section{Conclusion}\label{secCon}
Identification of large-scale linear dynamic systems may be subject to two relevant problems.
First, the number of unknown variables to be estimated may be large.
Using an MCMC scheme to reconstruct the posterior of the impulse response coefficients given the identification data, means that
any step of the algorithm can be computationally expensive (a matrix of very large dimension has to be inverted
to draw samples from the full conditional distributions of $\theta$).
Hence, it is important to derive schemes where small groups of variables
are updated. In addition, ill-conditioning can affect the problem due to input collinearity
which can also lead to slow mixing of the generated Markov chains.
This aspect requires a careful selection of the blocks to be updated.
An MCMC strategy based on the stable spline prior has been here proposed that addresses both of the above issues.
It relies on Gibbs sampling complemented with a strategy where overlapping blocks are updated.
The updating frequencies of such blocks are regulated by the level of collinearity among different system inputs.
Simulation results show good performance of the proposed algorithm.
In future work we plan to provide a theoretical analysis regarding the convergence of the proposed MCMC scheme
which will require to extend previous studies on Gibbs sampling like that reported in \cite{Robert95,RS97,convergnce14}.

\bibliographystyle{IEEEtrans}
\bibliography{biblio}

\end{document}